\documentclass[12pt,reqno]{amsart}
\usepackage{graphicx,times,bbm,url}
\usepackage[margin=1in]{geometry}

\usepackage{tikz}
\usetikzlibrary{arrows}
\tikzset{>=stealth'}

\usepackage[bookmarks=false,unicode,colorlinks,urlcolor=red!70!black,citecolor=red!70!black,linkcolor=blue]{hyperref}

\usepackage{subfig}
\newcommand{\RR}{\mathrm{I\kern-0.20emR}}
\newcommand{\D}{\mathrm{d}\kern0.2pt}

\newcommand{\vp}{{\varphi}}
\newcommand{\pd}[2][]{\frac{\partial #1}{\partial #2}}

\newcommand{\mref}[1]{%
\href{http://www.ams.org/mathscinet-getitem?mr=#1}{#1}}

\newcommand{\arxiv}[1]{%
\href{http://front.math.ucdavis.edu/#1}{ArXiv:#1}}

\title{Spilling from a cognac glass}
\author[T. Kulczycki]{Tadeusz Kulczycki}
\author[M. Kwa{\'s}nicki]{Mateusz Kwa{\'s}nicki}
\author[B. Siudeja]{Bart{\l}omiej Siudeja}
\address{T. Kulczycki and M. Kwa{\'s}nicki \\ Institute of Mathematics and Computer Science \\  Wroc{\l}aw University of Technology  \\ Wybrze{\.z}e Wyspia{\'n}skiego 27 \\ 50-370 Wroc{\l}aw, Poland.}
\address{B. Siudeja \\ Department of Mathematics \\ University of Oregon \\  Eugene, OR 97403. }
\email{Tadeusz.Kulczycki@pwr.wroc.pl}
\email{Mateusz.Kwasnicki@pwr.wroc.pl}
\email{siudeja@uoregon.edu}

\begin{document}

\maketitle

\section{Introduction}

\begin{figure}[b]
\centering
\begin{tikzpicture}
    \node[anchor=south west,inner sep=0] (image) at (0,0) {\includegraphics[width=0.5\textwidth]{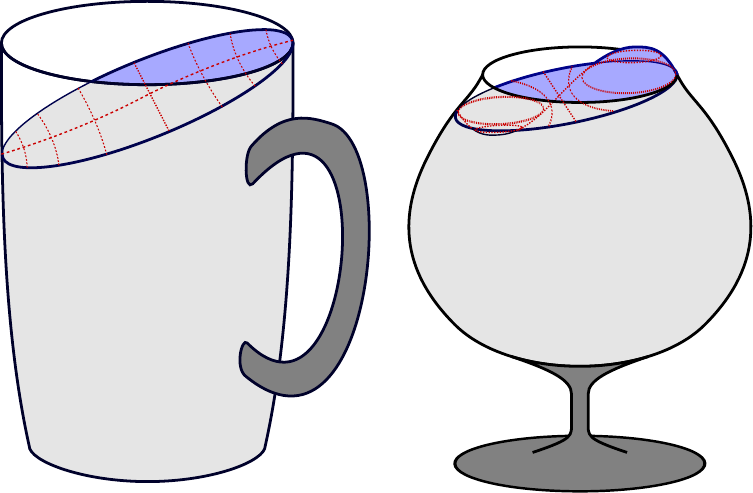}};
    \begin{scope}[x={(image.south east)},y={(image.north west)}]
      \draw (0.62,1) node (a) {{\footnotesize high spots}};
      \draw[->,red,very thick] (a.west) -- (0.393,0.92);
      \draw[->,red,very thick] (a.east) -- (0.85,0.91);
      %\draw (0.2,0) node [below] {\scriptsize (a)};
      %\draw (0.77,0) node [below] {\scriptsize (b)};
    \end{scope}
\end{tikzpicture}
    
\subfloat[A cup.\label{cupnglassa}]{\hspace{5cm}}
\subfloat[A snifter.\label{cupnglassb}]{\hspace{5cm}}

  \caption{High spots in a coffee cup and a snifter.}
  \label{cupnglass}
\end{figure}
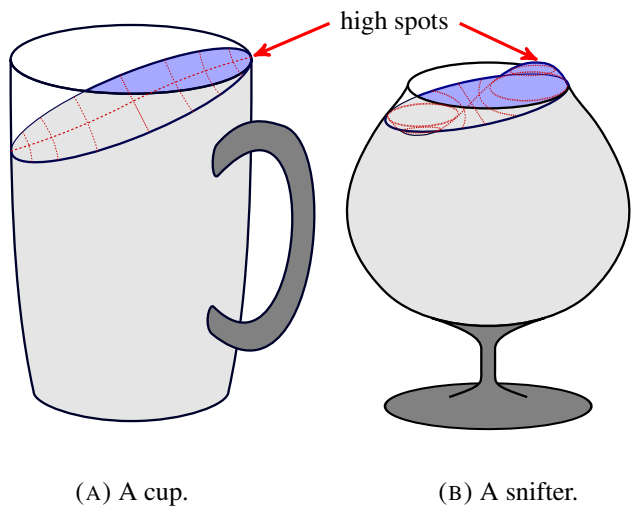

The 2012 Ig~Nobel Fluid Dynamics Prize was awarded to R.~Krechetnikov and H.~Mayer for their study of people walking while carrying a filled coffee mug \cite{KM2012}. They show that coffee spills so often because the sloshing mode with the lowest-frequency (most noticeable in practice) in a typical coffee mug tends to get excited during walking. Authors model oscillations of the coffee as appropriate mixed Steklov problem.

However, there is another reason for spilling from a mug: \emph{high spot} on the boundary. The maximal elevation of the lowest-frequency liquid oscillation in a typical coffee mug is located on the boundary (see \autoref{cupnglassa}). This effect, proved rigorously by Kulczycki and Kwa{\'s}nicki \cite{KK2012}, makes spilling even easier. On the other hand, in a bulbous snifter the lowest-frequency sloshing mode attains its maximal elevation (high spot) inside a snifter \cite{KK2012}, reducing the risk of spilling (see \autoref{cupnglassb}). 

The position of the high spot clearly depends on the container. Quite recently, this phenomenon has been studied by Faltinsen and Timokha \cite{FT2012}.  The natural limiting case of bulbous containers is an infinite ocean, covered with ice, with a single round hole. The corresponding sloshing problem is known as \emph{the ice-fishing problem}. Kozlov and Kuznetsov \cite{KK2004} showed that the maximal amplitude occurs approximately the third of the way from the boundary to the center of the hole, and the amplitude is over $50\%$ larger than at the boundary.

Note that sloshing dynamics is very important from the engineering point of view. Tanks of liquid propelled rockets include carefully designed baffles to mitigate the free surface effect. Those baffles minimize sloshing, preventing jitter (course deviations). Similar prevention techniques are used in ships, where water ballast, or liquid cargo can lead to capsizing. Even improper fire fighting technique can sink a ship. Sloshing is also used to minimize bouncing of the roller hockey ball. Interestingly, certain amounts of water in a ball can lead to resonance effects. 

In what follows we discuss only the high spot problem for sloshing. We begin with a formal definition of the model, present the explicit solutions for a cylindrical mug and plot some numerical approximations for harder cases. Then we describe a novel way of obtaining experimental evidence (see \autoref{photo}). We also present recently obtained rigorous results for solids of revolution. Finally we describe the relation of the high spot problem to the celebrated hot spots conjecture. 

\section{Mathematical model}

\begin{figure}[t]
\centering
\begin{tikzpicture}
    \node[anchor=south west,inner sep=0] (image) at (0,0) {\includegraphics[width=0.3\columnwidth]{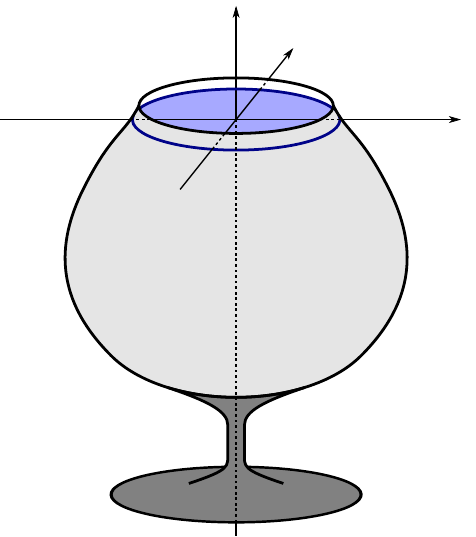}};
    \begin{scope}[x={(image.south east)},y={(image.north west)}]
    \draw (0.33020161,0.96190312/1.156) [below right=-2pt] node{\footnotesize $F$};
    \draw (0.85912897,0.4963678/1.156)node{\footnotesize $B$};
    \draw (0.62005341,0.53025341/1.156)node{\footnotesize $W$};
    \draw (0.95838757,0.9284346/1.156)node[above=-2pt]{\footnotesize $x$};
    \draw (0.6197406,1.06272564/1.156)node[above right=-3pt]{\footnotesize $y$};
    \draw (0.50797045,1.13862929/1.156)node[right=-2pt]{\footnotesize $z$};
    \end{scope}
\end{tikzpicture}

\caption{The domain for (\ref{slosh1}--\ref{slosh4}).}
  \label{notation}
\end{figure}

Small oscillations of ideal (that is, inviscid, incompressible and heavy) liquid with no surface tension are accurately described by the linear water-wave theory (see Fox and Kuttler \cite{FK83} and Troesch \cite{Tr60} for a historical overview, and \cite{GHLST,Ibr} for applications). In this framework, sloshing modes and frequencies correspond to solutions of the mixed Steklov eigenvalue problem (\cite{KK01,Mo64})
\begin{align}
 \Delta \vp & = 0  \text{ in $W$,} \label{slosh1} \\
 \pd[\vp]{z} & = \nu \vp  \text{ on $F$,} \label{slosh2} \\
 \frac{\partial \vp}{\partial \vec{n}} & = 0  \text{ on $B$,} \label{slosh3} \\
 \int_F \vp & = 0, \label{slosh4}
\end{align}
where $W$ is the part of the container filled with liquid (in its mean position), $B$ is the wetted part of its boundary, and $F$ is the free surface of the liquid (see \autoref{notation}). We choose a Cartesian coordinate system $(x,y,z)$ in which the $z$-axis is vertical and pointing upwards, and the free surface $F$ lies in the plane $z = 0$. The normal derivative on $B$ is denoted by $\frac{\partial}{\partial \vec{n}}$.

Any solution of the problem~(\ref{slosh1}--\ref{slosh4}) consists of an eigenfunction $\varphi$ and the corresponding eigenvalue $\nu$. In the language of hydrodynamics, there is a normal mode of the oscillations of the liquid such that the velocity of an element of liquid at point $(x,y,z)$ at time $t$ is equal to
\begin{equation}
\label{velocity}
  \cos(\omega t + \alpha) \nabla \varphi(x,y,z),
\end{equation}
where $\tfrac{1}{2 \pi} \omega = \tfrac{1}{2 \pi} \sqrt{\nu g}$ is the frequency of liquid oscillation, $g$ is the acceleration due to gravity and $\alpha$ is an arbitrary constant.

For sufficiently regular domains it is known that the mixed Steklov problem~(\ref{slosh1}--\ref{slosh4}) has a discrete sequence of eigenvalues
\[
  0 < \nu_1 \le \nu_2 \le \nu_3 \le \ldots \to \infty ,
\]
and the corresponding modes $\varphi_n \in H^1(W)$, $n= 1, 2, 3, \ldots$, restricted to the free surface $F$, form (together with a constant function) a complete orthogonal set in $L^2(F)$. In hydrodynamics, the first eigenfunction corresponding to the least eigenvalue $\nu_1$ plays an important role: it has the smallest decay rate due to non-ideal effects for real-life liquids.

The following observation plays a key role in our study: the amplitude of the oscillations of the liquid described by~(\ref{velocity}) on the free surface is proportional to $\tfrac{\partial \varphi}{\partial z}$, which, by~(\ref{slosh2}), is equal to $\nu \varphi$. Hence, the location of the high spot coincides with the location of the maximum of $|\varphi|$ on $F$, where $\varphi$ corresponds to $\nu_1$.

\section{Cylindrical mug --- a simple domain with explicit solution}\label{secmug}
 Krechetnikov and Mayer \cite{KM2012} modelled a mug by a cylinder 
\[
  W = \{(x,y,z): \, x^2 + y^2 < 1, \, z \in (-h,0)\} ,
\]
with the unit disk as free surface
\[
 F = \{(x,y,0): \, x^2 + y^2 < 1\}. 
\]
In this particularly nice case, the solutions of the mixed Steklov problem~(\ref{slosh1}--\ref{slosh4}) are known explicitly.
There are two linearly independent eigenfunctions corresponding to the least eigenvalue $\nu_1 = \nu_2$. In cylindrical coordinates $x = r \cos \theta$, $y = r \sin \theta$, these sloshing modes are given by
\begin{align}
  \varphi_1 & = J_1(j_{1,1}' r)  \cosh(j_{1,1}' (z  +  h)) \sin{\theta}, \label{eigf1}\\
  \varphi_2 & = J_1(j_{1,1}' r)  \cosh(j_{1,1}' (z  +  h)) \cos{\theta}, \label{eigf2}
\end{align}
where $J_1$ is the Bessel function of the first kind, and $j_{1,1}' \approx 1.8412$ is the first zero of $J_1'$. \autoref{fig7} shows a few level sets of  some linear combination of $\varphi_1$ and $\varphi_2$. Note that values on the top surface translate into liquid elevations.

Clearly $\varphi_1(x, y, 0)$ is an odd and increasing function of $y$, and attains extreme values (high spots) at boundary points $(0, 1)$ and $(0, -1)$. Similarly, $\varphi_2$ has its high spots at $(1, 0)$ and $(-1, 0)$. In fact, any linear combination of $\varphi_1$ and $\varphi_2$ will have high spots on the boundary, as on \autoref{fig7}.

\section{Numerical approximations}
\begin{figure}[t]
  \begin{center}
    \hspace{\fill}
    \subfloat[A cylinder.\label{fig7}]{
    \includegraphics[width=0.2\columnwidth]{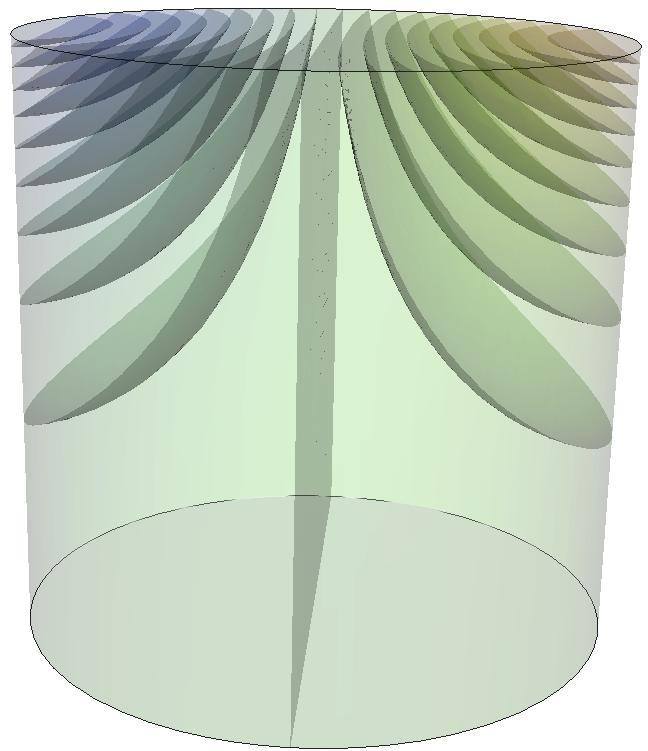}}
    \hspace{\fill}
    \subfloat[A silo.\label{figcos}]{
    \includegraphics[width=0.3\columnwidth]{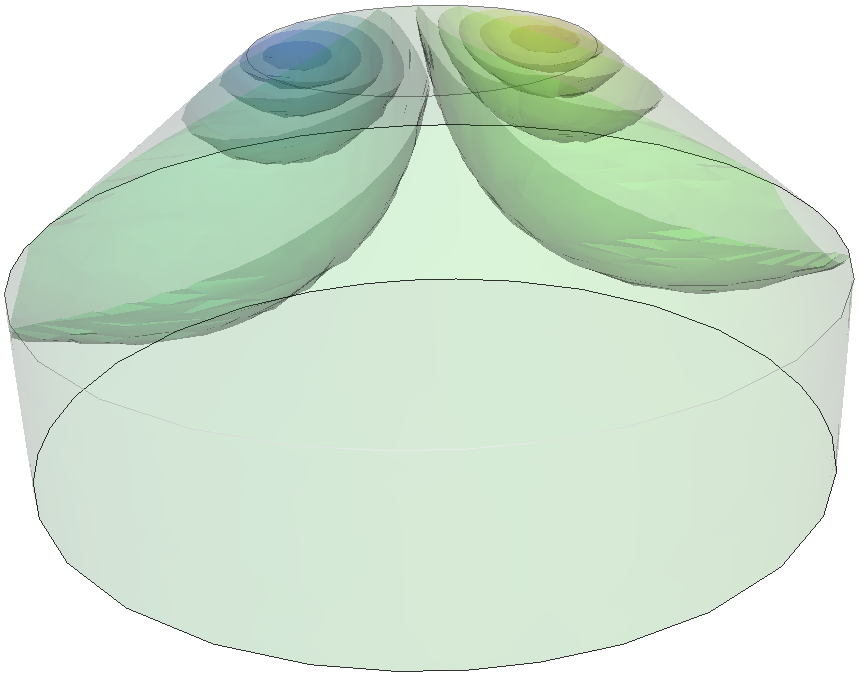}}
    \hspace{\fill}
    \subfloat[A trough.\label{figtr}]{
    \includegraphics[width=0.25\columnwidth]{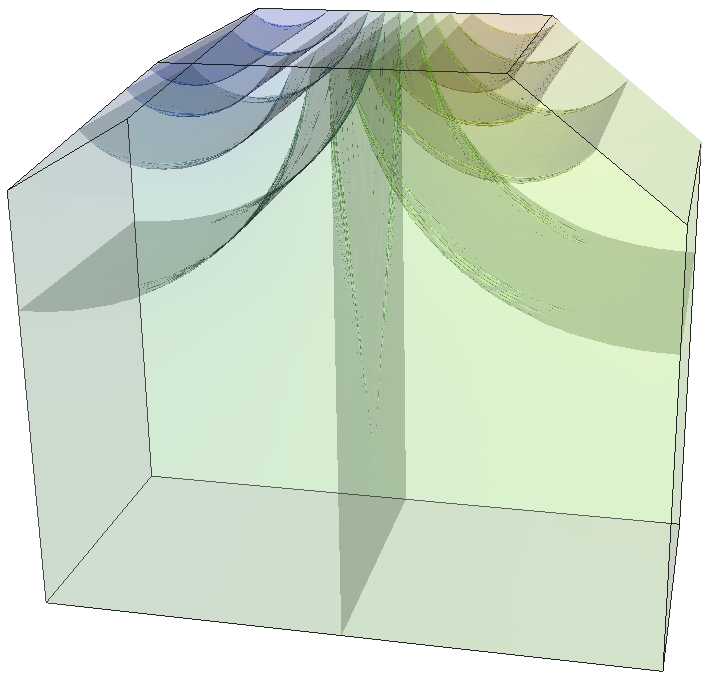}}
    \hspace{\fill}

  \end{center}
  \caption{Numerical solutions: level sets of fundamental eigenfunctions.}
\end{figure}

Finite Elements Methods can provide approximate positions for high spots for more complicated domains. Using FEniCS (\href{http://fenicsproject.org}{fenicsproject.org}, \cite{LMW}) we implemented a cylindrical mug and we obtained level sets clearly indicating that the extremal points are on the boundary (see \autoref{fig7}). Of course the eigenfunctions are explicit in this case, hence we could just plot the exact solutions from \autoref{secmug}. 

As a second example, we tried a grain silo shaped container, as on \autoref{figcos}. The high spot is clearly inside. This fact can actually be rigorously proved. See the discussion in \autoref{rigorous}.

Finally, we implemented a trough with a hexagonal cross-section and a small length (cf. \cite{KK2011}). \autoref{figtr} shows this trough with the level sets of the lowest-frequency mode. Clearly, the maximum is not on the boundary. 

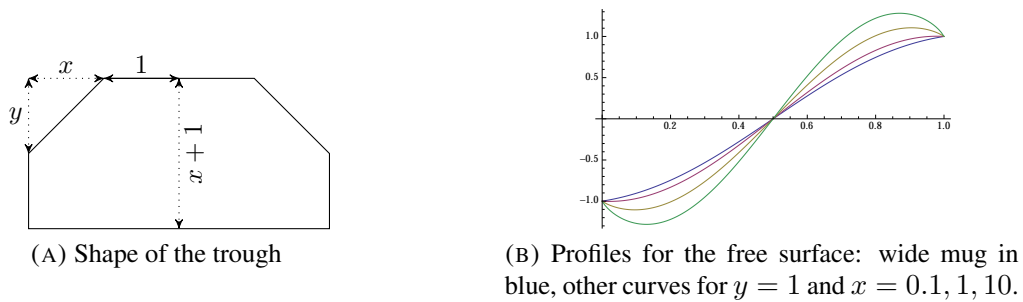
\begin{figure}[t]
  \begin{center}
    \hspace{\fill}
    \subfloat[Shape of the trough\label{shape}]{
    \makebox[0.45\textwidth]{
\begin{tikzpicture}[scale=1]
  \draw (0,0) -- (0,1) -- (1,2) -- (3,2) -- (4,1) -- (4,0) -- cycle;
  \draw[<->,dotted] (2,0) -- (2,2) node[below=-2pt,pos=0.5,sloped] {\footnotesize $x+1$};
  \draw[<->,dotted] (0,1) -- (0,2) node[left=-2pt,pos=0.5] {\footnotesize $y$};
  \draw[<->,dotted] (0,2) -- (1,2) node[above=-2pt,pos=0.5] {\footnotesize $x$};
  \draw[<->] (1,2) -- (2,2) node[above=-2pt,pos=0.5] {\footnotesize $1$};
\end{tikzpicture}}}
    \hspace{\fill}
\subfloat[Profiles for the free surface: wide mug in blue, other curves for {\small $y=1$} and {\small $x=0.1,1,10$}.\label{profiles}]{
\makebox[0.45\textwidth]{\includegraphics[width=0.3\columnwidth]{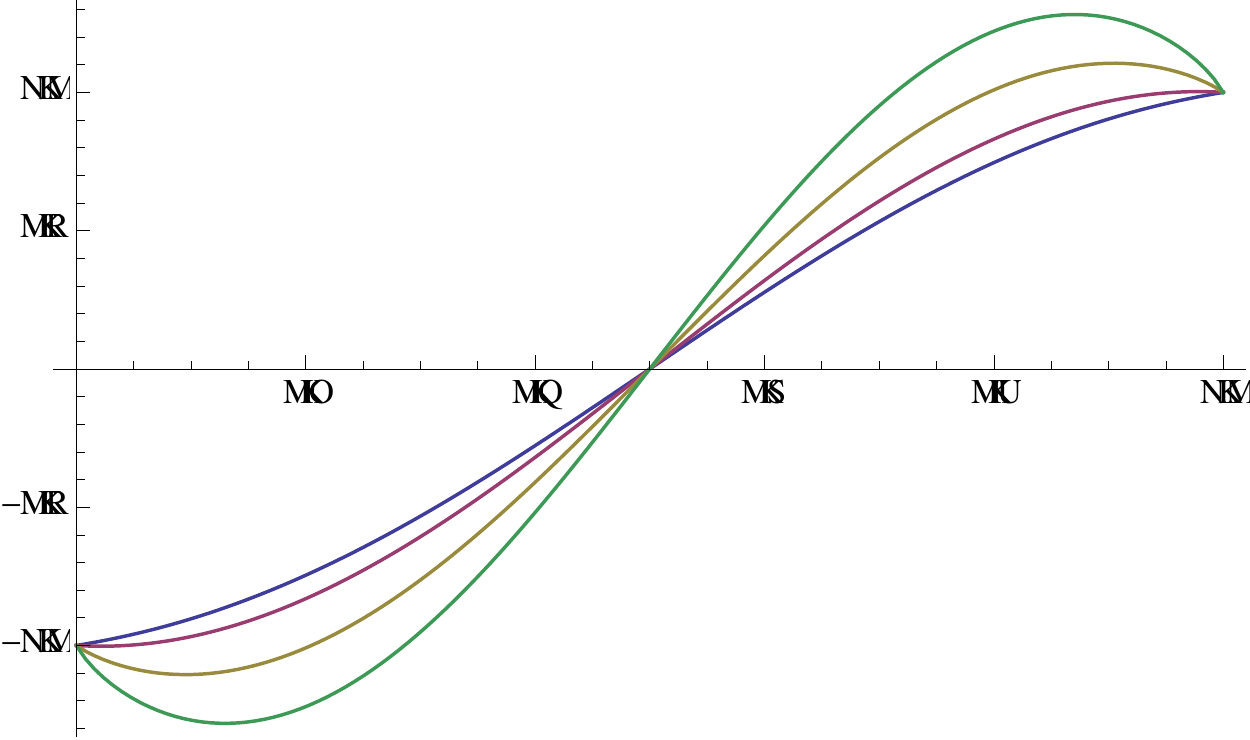}}}
    \hspace{\fill}

  \end{center}
  \caption{Sloshing in a trough}
\end{figure}

However, this trough and its lowest-frequency mode are essentially two dimensional. We exploited this reduction in dimension to obtain more accurate numerical profiles of the free surface shown on \autoref{profiles}. The blue curve corresponds to a trapezoidal trough that is 50\% wider at the top than at the bottom, and the maximum is on the boundary (as in coffee mug). Other profiles correspond to troughs from \autoref{shape} (cf. \autoref{figtr}) with different slopes. Generally, the maximum moves toward the center as the sloped part approaches horizontal line. The high spot also becomes more pronounced. Recently, similar numerical calculations have been made in \cite{FT2012} for spherical tanks (see Figure 2, ibid.).

\section{Physical experiment}\label{photo}

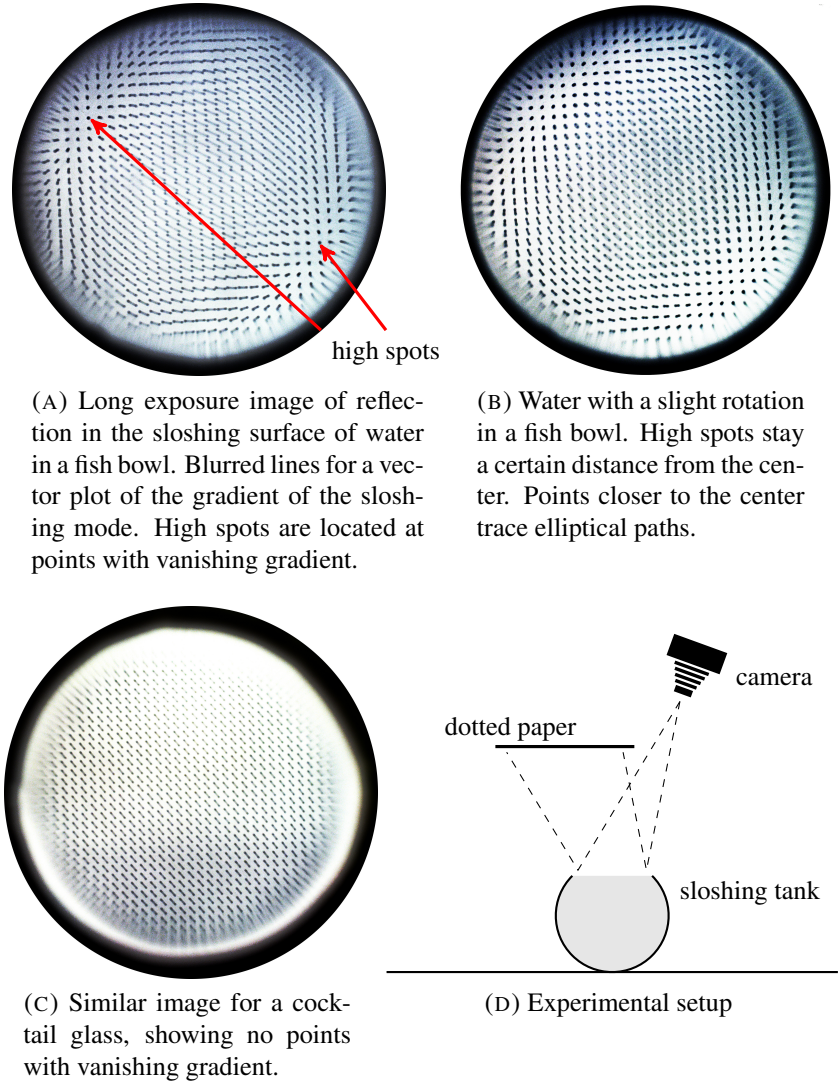
\begin{figure}[t]
  \centering
  \subfloat[Long exposure image of reflection in the sloshing surface of water in a fish bowl. Blurred lines for a vector plot of the gradient of the sloshing mode. High spots are located at points with vanishing gradient.\label{fishbowl}]{
\begin{tikzpicture}
    \node[anchor=south west,inner sep=0] (image) at (0,0) {\includegraphics[width=0.3\textwidth]{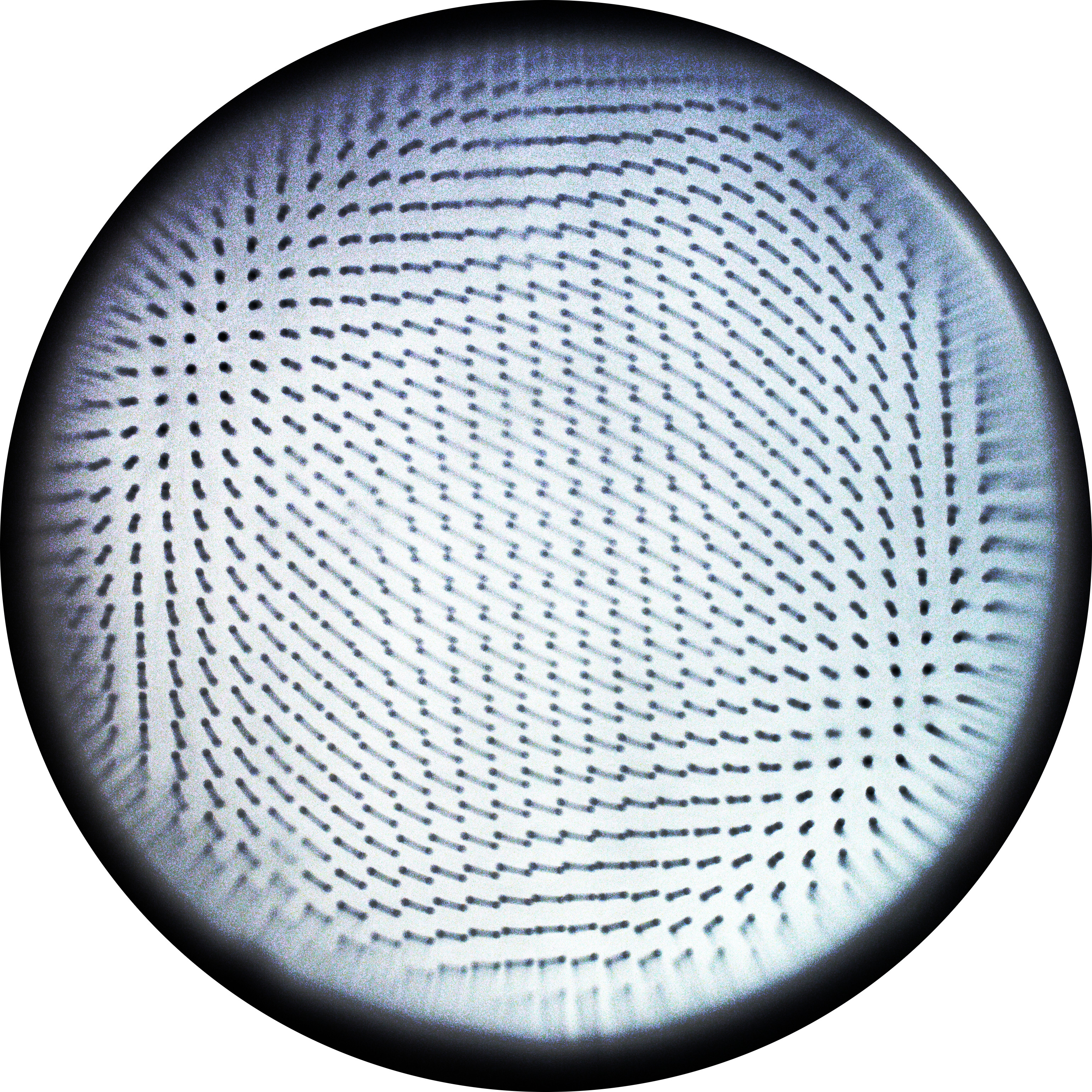}};
    \begin{scope}[x={(image.south east)},y={(image.north west)}]
      \draw (1,0) node [above] (a) {{\footnotesize high spots}};
      \draw[->,red,very thick] (a.north west) -- (0.21,0.688);
      \draw[->,red,very thick] (a.north) -- (0.824,0.352);
    \end{scope}
  \end{tikzpicture}}
  \subfloat[Water with a slight rotation in a fish bowl. High spots stay a certain distance from the center. Points closer to the center trace elliptical paths. \label{obrot}]{
  \includegraphics[width=0.3\textwidth]{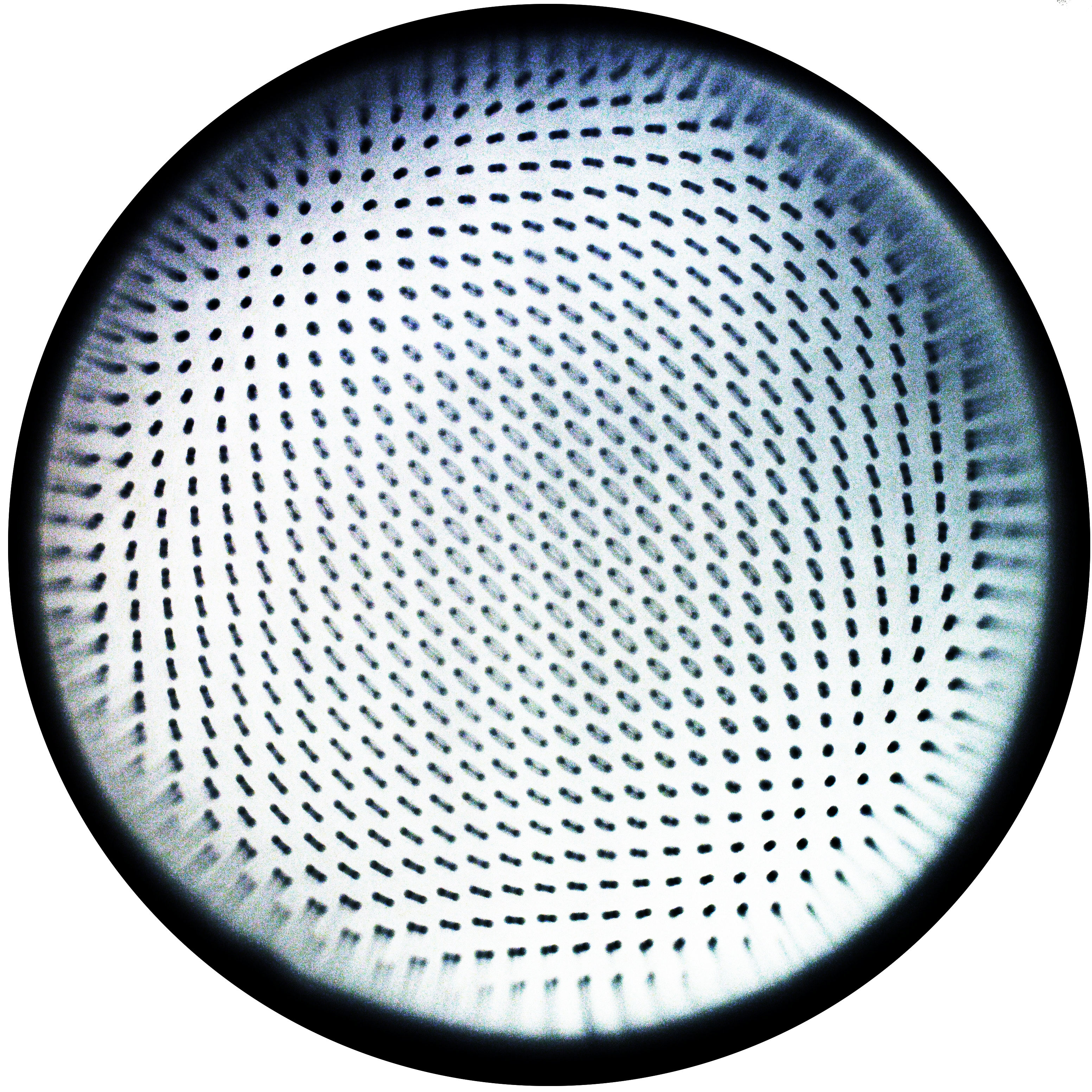}}

  \subfloat[Similar image for a cocktail glass, showing no points with vanishing gradient.\label{cocktail}]{
  \includegraphics[width=0.3\textwidth]{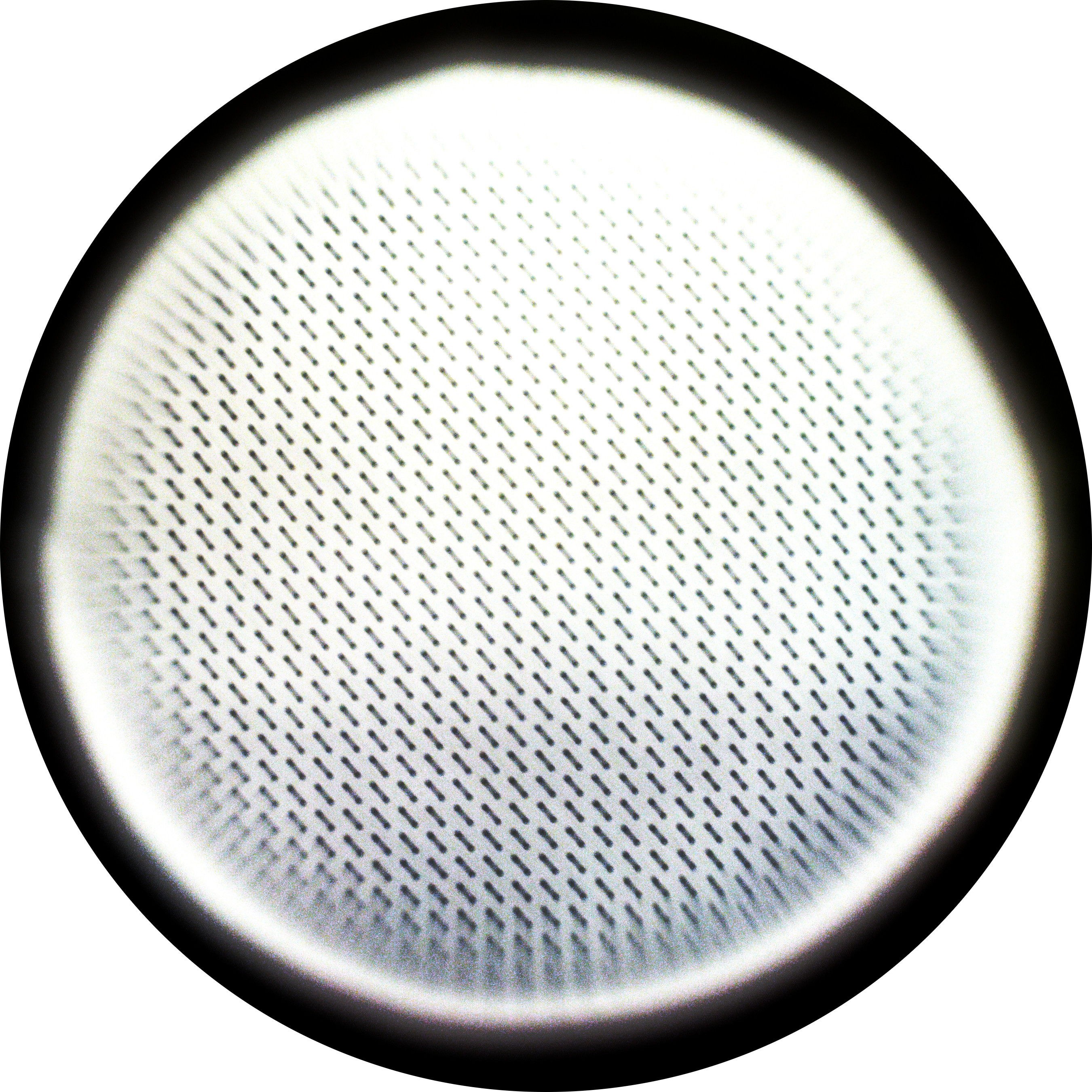}}
  \subfloat[Experimental setup\label{experiment}]{
\begin{tikzpicture}[scale=1.5]
    \draw[thick] (-2,0) -- (2,0); %floor
    \filldraw[fill=black!10!white,thick] (-0.35355,0.85355) arc (135:405:0.5); %tank
    \draw (0.5,0.7) node[right]{{\footnotesize sloshing tank}};
    \begin{scope}[shift={(0.55,3)}, rotate=250] %camera
      \fill (0,0) rectangle (0.2,0.5);
      \fill (0.22,0.09) rectangle (0.25,0.41);
      \fill (0.27,0.11) rectangle (0.30,0.39);
      \fill (0.32,0.13) rectangle (0.35,0.37);
      \fill (0.37,0.15) rectangle (0.40,0.35);
      \fill (0.42,0.17) rectangle (0.48,0.33);
    \end{scope}
    \draw (1,2.6) node[right]{{\footnotesize camera}};
    \draw[very thick] (-1.03539,2) -- (0.196246,2); %dotted paper
    \draw (-0.9,1.95) node[above]{{\footnotesize dotted paper}};
    \draw[dashed] (0.6,2.4) -- (0.30355,0.90355) -- (0.096246,1.95); %light
    \draw[dashed] (0.6,2.4) -- (-0.30355,0.90355) -- (-0.93539,1.95); %light
  \end{tikzpicture}}
  \caption{Photographic experiment}
\end{figure}

We have also performed experiments with bulbous and cylindrical containers. We tried to photograph the oscillations, and this proved difficult in any conventional way. Nonlinear effects caused by relatively large amplitude of oscillations and non-ideal nature of water are definitely noticeable. These problems are compounded by the  existence of two modes for the lowest frequency. This phenomenon causes whirling of the fluid if both modes are present and shifted in phase.

Instead, we photographed a reflection of a dotted piece of paper on a very slightly disturbed surface of the liquid (see \autoref{experiment}). This approach is consistent with the infinitesimal nature of the mathematical sloshing model. We would disturb the water, then wait until the surface became still to the naked eye. Then, the long exposure photograph would allow us to obtain mostly blurred images, with just a few clearly visible dots. Except at local extrema of the sloshing amplitude, planes tangent to the liquid surface oscillate, creating a path for each dot. On the other hand, at the extremal point the tangent plane is always horizontal, and the corresponding dot is sharp. Note that even this experiment is susceptible to mixing of the modes, hence whirling effect  (see \autoref{obrot}).

\autoref{fishbowl} shows the image obtained for a bulbous container (a fish bowl), while \autoref{cocktail} the image for a conical tank (a cocktail glass). There are clearly two extremal points away from the boundary in the bulbous container. On the other hand we obtained almost the same blurred paths for all dots in the conical container.

\section{Rigorous results}\label{rigorous}
\begin{figure}[t]
\centering
\subfloat[The domain\label{fig8}]{
\begin{tikzpicture}[scale=0.6]
  \draw[->] (-2.5,0) -- (2.5,0) node[above] {\footnotesize $x$};
  \draw[->] (0,-4.5) -- (0,1.5) node[right] {\footnotesize $z$};
  \draw[->] (-0.8,-1.1) -- (0.8,1.1) node[right] {\footnotesize $y$};
  \fill[fill=black!20!white, draw=black]
    (0,0) -- (0,-4) .. controls +(0.8,0) and +(-0.2,-2) .. (1.2,-0.64) -- cycle;
  \draw[thick] (-2,0) .. controls +(0.2, -2) and +(-1, 0) .. (0,-4) .. controls +(1, 0) and +(-0.2, -2) .. node[right] {\footnotesize \!$B$} (2,0);
  \draw[thick] (0,0) ellipse (2 and 0.8); 
  \node[above] at (-1,0.62) {\footnotesize $F$};
  \node at (-0.7,-2) {\footnotesize $W$};
  \end{tikzpicture}
}
\subfloat[Its profile\label{fig8b}]{
\begin{tikzpicture}[scale=0.6]
  \draw[->] (-2.5,0) -- (2.5,0) node[above] {\footnotesize $r$};
  \draw[->] (0,-4.5) -- (0,1.5) node[right] {\footnotesize $z$};
  \fill[fill=black!15!white,draw=black,thick]
    (0,0) -- (0,-4) parabola bend(0,-4) (2,0) -- (0,0);
  \draw (1,-0.1) node[above] {\footnotesize $F(D)$};
  \node at (0.8,-1.5) {\footnotesize $D$};
  \node at (1.95,-2.5) {\footnotesize $B(D)$};
\end{tikzpicture}
}
  \subfloat[A profile with symmetric eigenfunction\label{fig8c}]{
  \begin{tikzpicture}[scale=0.6]
  \draw[->] (-2.5,0) -- (5.5,0) node[above] {\footnotesize $r$};
  \draw[->] (0,-4.5) -- (0,1.5) node[right] {\footnotesize $z$};
  \fill[fill=black!15!white,draw=black,thick]
%   (0,-4) arc (-90:0:4) -- (3.5,0) arc (0:-85:3.5) -- ++(0,0.5) arc (-85:0:3) -| (0,-4); 
   (0,0) -- (0,-4) .. controls (2,-4) and (2,-0.1) .. (2.5,-0.1) .. controls (3,-0.1) and (3,-4) .. (4,-4) .. controls (5,-4) and (5,-2) .. (5,0);
  \end{tikzpicture}
  }
  \caption{Solids of revolution}
\end{figure}
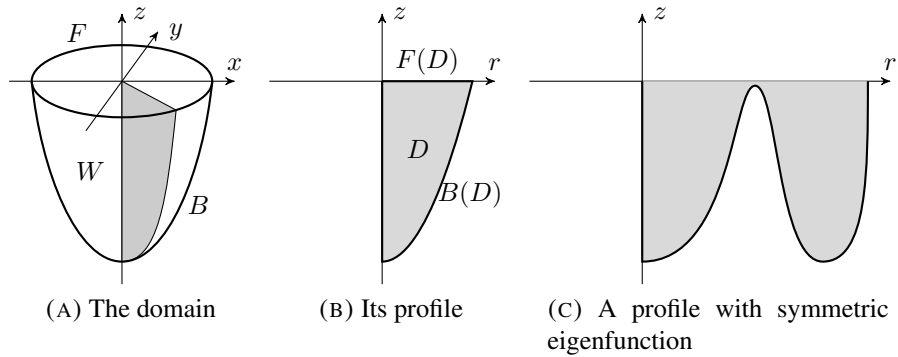

Finally, we discuss
rigorous results for solids of revolution. Let $W$ be the $3$-dimensional domain obtained by rotating a profile $D$ (\autoref{fig8b}) around the $z$-axis (\autoref{fig8}). Usually the free surface $F$ is a disk in the plane $z = 0$. Such containers were considered in  \cite{KK2012,GHLST,LBK}.

We assume that the modes corresponding to $\nu_1$ are antisymmetric. Many domains have this property. Nevertheless, strange examples of rotationally symmetric fundamental eigenfunctions exist, for example for a profile from \autoref{fig8c}.

Under the  antisymmetry  assumption, there are two linearly independent antisymmetric eigenfunctions corresponding to the least eigenvalue $\nu_1 = \nu_2$. In cylindrical coordinates they have the form
\begin{align*}
  \varphi_1 & = \psi(r,z) \cos \theta, & \varphi_2 & = \psi(r,z) \sin \theta,
\end{align*}
where $\psi(r,z)$ is a function defined on $D$.

By~\cite[Theorems~1.1 and 1.2]{KK2012}, if $W$ is a convex solid of revolution contained in the infinite cylinder $\{(x, y, z) : (x, y, 0) \in F\}$ (the last property is often called \emph{John's condition}), then indeed $\nu_1 = \nu_2$ correspond to antisymmetric eigenfunctions $\varphi_1$ and $\varphi_2$, and the high spots of the modes described by $\varphi_1$ and $\varphi_2$ are located on the boundary, e.g. for a profile from \autoref{fig8b}. On the other hand, \cite[Proposition~1.3]{KK2012} (see also \cite{FT2012}) asserts that if $B$ and $F$ form an obtuse angle then $\varphi_1(x,y,0)$, $\varphi_2(x,y,0)$ attain their maxima inside $F$, as shown on \autoref{cupnglassb} and \autoref{figcos}. 

The proof  of \cite[Theorem 1.2]{KK2012}  is based on the technique of domain deformation, and uses ideas of D.~Jerison and N.~Nadirashvili \cite{JN2000}. They studied related \emph{hot spots conjecture} (posed by J.~Rauch in 1974). Similarity between the high and hot spots is not an accident. It showcases the deep connection between sloshing problems and classical Laplace spectral problems with Neumann boundary conditions. We describe some of these in the next section.

\section{Relation to the classical spectral problems}
Traditionally, Laplace eigenvalues are intuitively understood via heat flow. The Dirichlet eigenvalues of a domain $D$ govern the heat on a plate $D$ with zero temperature on the boundary, while the Neumann eigenvalues give the insulated plate. This explains the name hot spots conjecture. The famous paper ``Can one hear the shape of drum?'' by Kac\cite{Kac} emphasizes another practical interpretation for Dirichlet spectrum. The corresponding eigenfunctions give the shapes of the drum membrane vibrating in one of the characteristic frequencies. The Neumann spectrum is harder to explain this way, as one would need a membrane that is free to move up and down, but is horizontal near the boundary. This behavior can however be attributed to the free surface of an ideal liquid in a container. Therefore, sloshing not only generalizes the Neumann spectral problem, but also provides physical intuition. In particular the hot spots problem might as well be called the high spot problem for a container with vertical walls.

Roughly speaking, the hot spots conjecture states that in a thermally insulated domain, for ``typical'' initial conditions, the hottest point will move towards the boundary of the domain as time passes. The mathematical formulation of the hot spots conjecture is the following: the extrema of every fundamental eigenfunction of the Neumann eigenvalue problem for a domain $D$ are located on the boundary of $D$. The conjecture was proved for sufficiently regular planar domains by Ba\~nuelos and Burdzy \cite{BB1999}, Jerison and Nadirashvili \cite{JN2000}, but disproved for some domains with holes by Burdzy and Werner \cite{BW99} and just one hole by Burdzy \cite{B2005}. See \emph{Nature} article~\cite{S1999} for more information, and Terence Tao's Polymath project~\cite{Tpoly} for current developments.

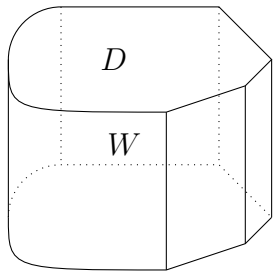
\begin{figure}
\begin{center}
\begin{tikzpicture}[scale=0.7]
  \draw (-1,0) arc (180:90:1) -- (3,1) -- (4,0) -- (3.5,-0.5) -- (2,-1) .. controls (-1,-1) .. (-1,0);
  \draw[yshift=-3cm, dotted] (-1,0) arc (180:90:1) -- (3,1) -- (4,0); -- (3.5,-0.5) -- (2,-1) .. controls (-1,-1) .. (-1,0);
  \draw[yshift=-3cm] (4,0) -- (3.5,-0.5) -- (2,-1) .. controls (-1,-1) .. (-1,0);
  \draw (-1,0) -- (-1,-3);% node [left,pos=0.5] {$L$};
  \draw[dotted] (0,1) -- (0,-2);
  \draw[dotted] (3,1) -- (3,-2);
  \draw (4,0) -- (4,-3);
  \draw (3.5,-0.5) -- (3.5,-3.5);
  \draw (2,-1) -- (2,-4);
  \draw (1.2,-1.2) node [below] {$W$};
  \draw (1,0) node {$D$};
\end{tikzpicture}
\end{center}
\caption{A sloshing cylinder $W$ with a cross-section $D$.
}
\label{cylindrical}
\end{figure}

The function $\varphi$ is a sloshing mode in the cylinder
\[
 W = \{ (x, y, z) : (x, y) \in D, \, z \in (-h, 0) \}
\]
if and only if $\psi(x, y) = \varphi(x, y, 0)$ is a Neumann eigenfunction on $D$ (see \autoref{cylindrical}). Furthermore
\[
  \varphi(x, y, z) = \psi(x, y) \cosh(\sqrt{\mu} (z + h))
\]
and $\nu = \sqrt{\mu} \tanh(\sqrt{\mu}h)$, where $\mu$ is the Neumann eigenvalue corresponding to $\psi$ and $\nu$ is the sloshing eigenvalue corresponding to $\varphi$ by \cite[Proposition 3.1]{KK2009}. This directly embeds the Neumann problem in the sloshing problem. Due to the very explicit connection many results and conjectures for Neumann spectrum can be quickly translated to the sloshing language. In particular, the high spots problem~(\ref{slosh1}--\ref{slosh4}) generalizes the hot spots conjecture. Note that we already discussed a special case of a mug in \autoref{secmug}. In (\ref{eigf1}-\ref{eigf2}) one may recognize Neumann eigenfunctions of the disk, given by the product of a Bessel and a trigonometric function. 

By the above identification and the results of~\cite{BB1999, JN2000}, for any convex domain $D$ with two orthogonal axes of symmetry (e.g. ellipses), the high spots of the fundamental sloshing modes for cylindrical tanks with cross-section $D$ are located on the boundary. Note that the counterexample for the hot spots conjecture given by Burdzy and Werner \cite{BW99} is rather similar to the example from \autoref{fig8c}. Indeed, their domain consists of a disk and a larger annulus, connected using spokes. The profile from \autoref{fig8c} shows two deep containers with a disk and an annulus as free surfaces, connected using extremely shallow portion. In both cases one would like to have two disjoint subdomains, even though that is not allowed. In the high spot problem there are more ways to make the domain connected due to one extra dimension. For that reason it is much easier to construct intuitive counterexamples for high spots problem.

The connection between Neumann problem and sloshing was also used to conjecture bounds for the eigenvalues. In particular the Dirichlet to Neumann eigenvalue comparisons \cite{LW86,Fr95,Fi04} was already extended to mixed Steklov eigenvalues \cite{BKPS10}. Furthermore, bounds for spectral functionals for Neumann eigenvalues naturally extend to the cylindrical sloshing case (see \cite[Section 13]{LS}), opening a whole new area of research for general containers. Finally, sloshing problem can be used to tackle, seemingly unrelated, spectral problems for nonlocal fractional Laplacian  \cite{BK04,KKMS2010, K2012}.

\end{document}